\documentclass[a4paper,10pt]{article}
\usepackage{amsfonts,amsmath,amsthm,amssymb}
\usepackage[colorlinks]{hyperref}

\newtheorem{thm}{Theorem}
\newtheorem{lem}[thm]{Lemma}

\theoremstyle{definition}
\newtheorem*{defn}{Definition}

\DeclareMathOperator{\diam}{{\operatorname{diam}}}
\DeclareMathOperator{\proj}{\operatorname{proj}}
\newcommand{\eps}{\varepsilon}
\newcommand{\R}{{\mathbb R}}
\newcommand{\Q}{{\mathbb Q}}

\title{TUBE-MEASURABILITY}
\author{Marianna Cs\"ornyei \and Laura Wisewell\thanks{Supported by a Royal Society Dorothy Hodgkin Fellowship}}
\date{}

\begin{document}
\maketitle
In this note we investigate the measurable sets of an outer measure defined as follows.
\begin{defn}
In $\R^n$, let $T_i$ denote an infinite tube of cross-sectional radius $r_i>0$; that is, the closed $r_i$-neighbourhood of some straight line.  For a set $E\subseteq\R^n$ we define its \emph{tube-measure} by
\[
  \mu(E):=\inf\left\{\sum_i \gamma_{n-1}r_i^{n-1}:\bigcup_i T_i\supseteq E\right\},
\]
where $\gamma_{n-1}$ is the volume of the unit ball of $\R^{n-1}$.
Call $E$ \emph{tube-null} if $\mu(E)=0$.
\end{defn}
A closely-related outer measure has been introduced by Carbery, Soria and Vargas in connection with Fourier localisation: they showed that every tube-null set is a ``set of divergence'' for the localisation problem (see \cite{csv:localisation}). They observed that the tube-measure is very badly-behaved, in that Borel sets need not be measurable.  Our main result is the following.
\begin{thm}
In all dimensions, the only tube-measurable sets are the tube-null sets and their complements.
\end{thm}
For the proof we will need some estimates for the tube-measure of sets.  Exact values are not known even for simple sets such as balls, except in the case $n=2$ which corresponds to the famous Plank Problem \cite{bang:plank}.
\begin{lem}
\label{lem:bounds}
For every set $E\subseteq\R^n$ we have the upper bound
\[
  \mu(E) \leq \min|\proj(E)|,
\]
where $\proj(E)$ denotes a projection of $E$ onto an $(n-1)$-dimensional subspace, and $|\,\cdot\,|$ denotes the Lebesgue outer measure in $\R^{n-1}$.
For bounded sets $E$ we also have the lower bound
\[
  \mu(E) \geq \frac{|E|}{\diam(E)},
\]
where this time $|\,\cdot\,|$ denotes Lebesgue outer measure in $\R^n$.
\end{lem}
The plank problem tells us that for convex $E\subseteq\R^2$ we actually have $\mu(E) = \min|\proj(E)|$.  The proof of our theorem could be much simplified if we could assert this in higher dimensions, even just for balls.  However the statement is not true for all convex bodies, with the tetrahedron in $\R^3$ providing a counterexample \cite{bang:plank}.
\begin{proof}[Proof of Lemma~\ref{lem:bounds}]
The upper bound is obvious, by covering $E$ with parallel tubes.
For the lower bound, note that for any tube $T$ of cross-sectional radius $r$ we have $|E\cap T|\leq\diam(E)\gamma_{n-1}r^{n-1}$.  So if $\bigcup_{i=0}^\infty T_i\supseteq E$, we have
\[
  |E|
  =\left|E\cap\bigcup_{i=0}^\infty T_i\right|
  \leq\sum_{i=0}^\infty|E\cap T_i|
  \leq\diam(E)\sum_{i=0}^\infty \gamma_{n-1} r_i^{n-1},
\]
from which the inequality follows by taking the infimum.
\end{proof}
Next we show that the upper bound is the exact value in the case of sets that are Cartesian products with $\R$.
\begin{lem}\label{lem:product}
If $A\subseteq\R^{n-1}$, then $\mu(A\times\R)=|A|$.
\end{lem}
\begin{proof}
The upper bound is immediate from the previous lemma.  For the lower bound, observe that for all $R>0$ we have
\[
\mu(A\times\R)\geq\mu(A\times[-R,R])\geq\frac{2R|A|}{2R+\diam(A)}
\]
by the previous lemma, which goes to $|A|$ as $R\to\infty$.
\end{proof}
In particular, the $\mu$-measure of a single tube is exactly its cross-sectional area, as expected.  We are now ready to prove the theorem.

\begin{proof}[Proof of the main theorem]
Let $E\subseteq\R^n$ be $\mu$-measurable and, for a contradiction, suppose that both $\mu(E)>0$ and $\mu(\R^n\setminus E)>0$.  Choose a ball for which $\mu(E\cap\text{ball})>0$.  Then for any $\eps>0$ we can find a family of tubes $T_i$ covering $E\cap\text{ball}$ such that
\[
  (1-\eps)\sum_i\mu(T_i)<\mu(E\cap\text{ball})\leq\sum_i\mu(E\cap\text{ball}\cap T_i).
\]
Therefore there exists a tube $T=T_i$ with
\[
  (1-\eps)\mu(T)\leq\mu(E\cap\text{ball}\cap T)\leq\mu(E\cap T).
\]
Subdivide this $T$ into the union of countably many non-overlapping ``square tubes'' $T=\bigcup R_i$ (where each $R_i$ is a shifted and rotated copy of $[-\delta_i,\delta_i]^{n-1}\times\R$). Without loss of generality we can assume that $\delta_i\in \Q$ for all $i$.  Then, using
Lemma~\ref{lem:product},
\begin{equation}\label{R}
(1-\eps)\mu(T)=(1-\eps)\sum_i\mu(R_i),
\end{equation}
hence
\[
  (1-\eps)\sum_i\mu(R_i)\le\mu(E\cap T)\le\sum_i\mu(E\cap R_i).
\]
Therefore we can choose a square tube $R=R_i$ with
\begin{equation}\label{RR}
  (1-\eps)\mu(R)\leq\mu(E\cap R).
\end{equation}
We can similarly do this for $\R^n\setminus E$.  The two square tubes we have found may be of different widths, but since both the widths are rational we can subdivide as before into the union of square tubes of some common smaller width $\delta>0$, and from each of the two collections select a tube  still satisfying \eqref{RR}.

So we now have some $\delta>0$ (which depends on $\eps$) and two copies of
$[-\delta,\delta]^{n-1}\times\R$ that we denote by $R_1$ and $R_2$, such that
\begin{equation}\label{twotubes}
\begin{aligned}
(1-\eps)\mu(R_1)&\leq\mu(R_1\cap E)\\
(1-\eps)\mu(R_2)&\leq\mu(R_2\setminus E).
\end{aligned}
\end{equation}
Choose disjoint $\delta$-balls $B_1\subseteq R_1$, $B_2\subseteq R_2$. We now want to pass from balls to very eccentric sets, because for such sets the upper and lower bounds of Lemma~\ref{lem:bounds} are almost equal. So into each ball $B_i$ ($i=1,2$) place a cuboid $C_i$ of diameter $2\delta$ with $n-1$ of its edges all equal to some small $\eta>0$ to be chosen later. By Pythagoras' Theorem, their measure is
\[
  |C_1|=|C_2|=\eta^{n-1}\sqrt{4\delta^2-(n-1)\eta^2}.
\]
Orient the two cuboids within the balls so that they both lie in the Cartesian product of $\R$ with a common cube of side $\eta$.  Then
\begin{align*}
\eta^{n-1}&\geq
\mu(C_1\cup C_2) &\text{by Lemma~\ref{lem:bounds}}\\
&= \mu\big((C_1\cup C_2)\cap E\big)+\mu\big((C_1\cup C_2)\setminus E\big) &\text{by measurability of $E$}\\
&\geq \mu(C_1\cap E)+\mu(C_2\setminus E)&\text{by monotonicity}\\
&= \mu(C_1)-\mu(C_1\setminus E)+\mu(C_2)-\mu(C_2\cap E)&\text{by measurability of $E$}.
\end{align*}
Now by Lemma~\ref{lem:bounds} we have $\mu(C_2)=\mu(C_1)\geq |C_1|/\diam(C_1)=|C_1|/2\delta$.  Also by measurability of $E$ and \eqref{twotubes}
\begin{gather*}
\mu(C_1\setminus E)\leq\mu(R_1\setminus E)\leq \eps\mu(R_1)=\eps(2\delta)^{n-1}\\
\mu(C_2\cap E)\leq\mu(R_2\cap E)\leq \eps\mu(R_2)=\eps(2\delta)^{n-1}.
\end{gather*}
Putting these together we find
\[
  \eta^{n-1}\geq 2\left(\frac{|C_1|}{2\delta}-\eps(2\delta)^{n-1}\right),
\]
that is,
\begin{equation}\label{2}
1\geq 2\left(\frac{|C_1|}{2\delta\eta^{n-1}}-\eps(2\delta/\eta)^{n-1}\right).
\end{equation}
We must choose suitable values of the parameters so that this is a contradiction. Note that
\[
  |C_1|/2\delta\eta^{n-1}=\sqrt{1-(n-1)(\eta/2\delta)^2}.
\]
So we choose our $\eps$ and $\delta$ as follows: first let
\[
  p<\sqrt{\frac{3}{4(n-1)}}
\]
so that
\[
  \sqrt{1-(n-1)p^2}>1/2.
\]
Choose $\eps$ so small that
\[
  \sqrt{1-(n-1)p^2}-\eps/p^{n-1}>1/2.
\]
Then in the above construction use this $\eps$ and let $\eta=2\delta p$. Then \eqref{2} gives
\[
1\geq2\left(\sqrt{1-(n-1)(\eta/2\delta)^2}-\eps(2\delta/\eta)^{n-1}\right)>1,
\]
a contradiction.
\end{proof}

\bigskip

\begin{tabular}{l@{\qquad\qquad}l}
Marianna Cs\"ornyei   & Laura Wisewell\\
Department of Mathematics  & Department of Mathematics\\
University College London & University of Glasgow\\
Gower Street & University Gardens\\
London WC1E 6BT & Glasgow G12 8QW\\
UK & UK\\
(\url{mari@math.ucl.ac.uk})&(\url{lw@maths.gla.ac.uk})
\end{tabular}
\end{document}